%% file: main.tex
\documentclass[letterpaper, 10 pt, conference]{ieeeconf}  

\IEEEoverridecommandlockouts                              
\usepackage{graphicx,amstext,amsmath,amssymb,color,psfrag,mathabx,cite,booktabs,mathtools,enumerate}
\usepackage{url}            

\definecolor{blue}{rgb}{0.0, 0.0, 1.0}

\newcommand{\R}{\mathbb{R}} 
\newcommand{\N}{\mathbb{N}} 

\newtheorem{thm}{Theorem}
\newtheorem{lem}{Lemma}
\newtheorem{ass}{Assumption}

\usepackage{tikz}
\usetikzlibrary{calc,positioning,decorations.pathreplacing}

\title{\LARGE \bf
Adaptive control mechanisms in gradient descent algorithms
}

\author{Andrea Iannelli
\thanks{The author is with the University of Stuttgart, Institute for Systems Theory and Automatic Control, 70550 Stuttgart, Germany
        {\tt\small iannelli@ist.uni-stuttgart.de.}}%
}

\begin{document}

\maketitle
\thispagestyle{empty}
\pagestyle{empty}

\begin{abstract}

The problem of designing adaptive stepsize sequences for the gradient descent 
method applied to convex and locally smooth functions is studied. We take an adaptive control perspective 
and design update rules for the stepsize that make use of both past (measured) and future (predicted) 
information. We show that Lyapunov analysis can guide in the systematic design of adaptive parameters striking a balance 
between convergence rates and robustness to computational errors or inexact gradient information. Theoretical and numerical results 
indicate that closed-loop adaptation guided by system theory is a promising approach for designing new classes of adaptive optimization algorithms 
with improved convergence properties.

\end{abstract}

\section{INTRODUCTION}

Convex optimization algorithms are at the core of many established methodologies in control, reinforcement learning and machine learning, for example receding horizon control \cite{Borrelli2017}, convex Q-learning \cite{Lu_CDC23_CvxQ}, and some classes of regularized neural networks training \cite{Ergen_Allerton_2019}. In all these applications, one typically requires algorithms that are fast (to reduce computational time) and robust (e.g., to be less sensitive to error in the problems data). Developing systematic methods to analyze these properties and design for them is a very important step to make these tools of more widespread and dependable use in applications. Because iterative optimization algorithms can be seen as open dynamical systems, tools and viewpoints from control theory can be a valid standpoint to approach these tasks \cite{Scherer_CDC23,SystemsTheoryofAlgorithms}.
Towards this goal, we consider here the basic unconstrained optimization setting where the objective $f$ is convex and only locally smooth. We approach the design of a stepsize sequence $(\alpha_k)_{k \in \N}$ as an adaptive control problem where the goal is to guarantee that the interconnection between the algorithm and the stepsize law converges to the minimizer set. We do this while capturing performance and robustness trade-off of this adaptive closed-loop.

Even though gradient descent (GD) methods  are standard in optimization, analysis and design of varying stepsizes is an active area of research. For the case of $L_\text{s}$-smooth objective, one line of work has developed pre-defined sequences of \emph{large} stepsizes (i.e., with instances where these are larger than $\frac{2}{L_\text{s}}$) and showed that they can accelerate convergence \cite{Altschuler2024,Grimmer_SJO}. While convergence guarantees and the sequence of stepsizes generally depend on the value of the stopping time which must be selected a-priori, very recently \cite{Zhang-arxiv24} showed that this strategy provably achieves anytime convergence guarantees that strictly improves upon the classic $\mathcal{O}(\frac{1}{k})$. 
Besides the restriction to smooth objectives, the fundamental idea of these approaches is to pre-compute the sequence of stepsizes independently of $f$ (except for its smoothness constant) and of the initial iterate, thus the stepsize sequence is effectively an open-loop input to the GD dynamics. In another line of work, the case of locally smooth objective has been addressed by proposing stepsizes that adapt to the local geometry \cite{Latafat2024_MathProg,latafat24aL4DC,Malitsky2024_AProx}, including extension to the proximal gradient method for composite problems. In these works, feedback from the current and past iterates is leveraged to estimate a sequence of local smoothness constants which are used in the stepsize's update. By doing so, the standard requirement of global smoothness is lifted, and progress to the optimal value is empirically much faster even though a better rate than $\mathcal{O}(\frac{1}{k})$ is not proved.   

Inspired by this prior work, we make a first attempt to combine feedback and feedforward (or predictive) actions in the selection of stepsize laws for gradient-based methods applied to locally smooth objective. We approach the problem similarly to the 
analysis-informed design of adaptive controllers whereby appropriately constructed Lyapunov functions guide the selection of parameters so that boundedness and convergence guarantees can be established. We also investigate questions on robustness of such adaptive systems, and recognize components of the designed adaptive mechanisms that can mediate between performance and robustness, e.g., to errors in the gradient information.
Besides the new technical results, we see as central contribution of this work the consideration of the role of feedback and predictive actions in adaptive optimization algorithms. Even though in a different context, connections between convex optimization and adaptive control were also studied in \cite{raginsky2010online}.

The paper is structured as follows. We discuss preliminaries on analysis of GD with varying stepsizes in Section \ref{sec:Preliminaries}.
In Section \ref{sec:Adapting}, we introduce our proposed algorithm and analyze and discuss its properties. We use a logistic regression example to show various features of the algorithm and compare it with other non-adaptive and adaptive algorithms in Section \ref{sec:Numerics}.
The proofs of all the new technical results are presented in the Appendix for the sake of readability. 

\emph{Notation:} 
We denote by $\langle x,y\rangle$ the standard Euclidean inner product in $\R^n$ and by $\| \cdot \|$ its induced norm. We use $\N$ for the set of natural numbers.
Given $u,v \in \R^n$, we refer to the following as the Pythagoras identity
\begin{equation*}
\|u\|^2=\|u-v\|^2-\|v\|^2+2\langle u,v\rangle
\end{equation*}
The convex hull of a set of points is denoted by $\text{conv}$.
\newpage
\section{Preliminaries}\label{sec:Preliminaries}
We consider the unconstrained problem
\begin{equation}\label{Problem_start}
	\min_{x \in \R^n} \quad\hspace*{-0.5em} f(x)
\end{equation}
and we denote with $X^\star$ its set of minimizers and with $f^\star$ its optimal value. 
\begin{ass}\label{Standing_ass}
We make the following standing assumptions:
\begin{itemize}

\item $f:\R^n\rightarrow \R$ is \emph{locally smooth}, i.e., $f$ is differentiable and, for every convex and compact set $D \subset \R^n$, there exists $L_D \in(0,\infty)$ such that $\forall x,y \in D$:
\begin{equation}\label{f_loc_smooth}
\|\nabla f(y)-\nabla f(x)\|\leq L_D \|y-x\|.  
\end{equation}

\item $f$ is \emph{convex}, i.e., $\forall x,y \in \R^n$
\begin{equation}\label{f_convex}
f(y) \geq f(x)+\langle\nabla f(x),y-x\rangle. 
\end{equation}

\item $X^\star  \neq \emptyset$ and $f^\star > - \infty$
\end{itemize}
\end{ass}
Local smoothness, or equivalently local Lipschitz continuity of the gradient, defines a rather general class of functions. For example, any twice differentiable $f$ is locally smooth since (\ref{f_loc_smooth}) holds with $L_D=\max_{v \in D} \| \nabla^2 f(v) \| $, which is finite due to continuity of the Hessian and compactness of $D$.

One of the most popular methods to solve (\ref{Problem_start}) is gradient descent (GD), which generates a sequence of iterate $(x_k)_{k \in \N}$ by applying the following simple recursion
\begin{equation}\label{GD_method}
x_{k+1}=x_k-\alpha_k \nabla f(x_k),\quad k \in \N,
\end{equation}
starting from a given $x_0$ and choosing an appropriate sequence of stepsizes \footnote{With a slight abuse of terminology, we will refer to \emph{time-varying} step-sizes whenever $\alpha$ changes across iterations, owing to the interpretation of (\ref{GD_method}) as a discrete-time dynamical system.} $(\alpha_k)_{k \in \N}$. We will denote by $F_k\coloneqq f(x_k)-f^\star$ the optimality gap at iteration $k$.

Standard convergence analyses of GD assume that $f$ is (globally) $L_\text{s}$-\emph{smooth}, i.e., $f$ is differentiable and there exists $L_\text{s} \in(0,\infty)$ such that $\forall x,y \in \R^n$
\begin{equation}\label{f_glob_smooth}
\|\nabla f(y)-\nabla f(x)\|\leq L_\text{s} \|y-x\|, 
\end{equation}
or equivalently $\forall x,y \in \R^n$
\begin{equation}\label{f_glob_smooth_2}
f(y) \leq  f(x)+\langle\nabla f(x),y-x\rangle + \frac{L_\text{s}}{2}\|y-x\|^2. 
\end{equation}
In this case, one can guarantee global convergence of (\ref{GD_method}) to an element of $X^\star$ by restricting the choice of stepsize \cite{Nesterov_Intro_2014}. Most of the analyses show convergence with constant stepsizes in the ranges $L_\text{s} \alpha \in(0,1]$ or $L_\text{s} \alpha \in[1,2)$. As a summary of the available analysis results, we provide a Lyapunov-based analysis that encompasses any time-varying stepsize sequence satisfying\footnote{We make the last two technical restrictions to simplify some steps in the derivation of the rates in view of the very unrestricted range of time-varying stepsizes; this is without loss of generality, e.g., \cite[Section 4]{Teboulle2022} consider the specific case where $L_\text{s} \alpha_k \in [1,2)$, $ \lim_{k \to +\infty} L_\text{s} \alpha_k \rightarrow 2$.}:  $L_\text{s} \alpha_k \in (0,2)$; $ \limsup_{k \to +\infty} L_\text{s} \alpha_k \neq 2$; $ \liminf_{k \to +\infty} \alpha_k \neq 0$. 

\begin{thm}\label{thm:GD_LF_smooth} 
Let $(x_k)_{k \in \N}$ be a sequence generated by the GD method (\ref{GD_method}) applied to a $L_\text{s}$-smooth function $f$ also satisfying Assumption \ref{Standing_ass} with any time-varying stepsize sequence satisfying $L_\text{s} \alpha_k \in (0,2), k \in \N$. Define the function
\begin{equation}\label{LF_smooth}
V^{\text{s}}_k(x^\star)\coloneqq\|x_k-x^\star\|^2,\quad x^\star \in X^\star.
\end{equation}
Then for any $x^\star \in X^\star$, $k \in \N$ and $x_0 \in \R^n$, it holds\footnote{We omit the argument of $V^{\text{s}}$ for brevity and formatting reasons.}
\begin{subequations}\label{GD_LF_smooth}
\begin{align}
V^{\text{s}}_{k+1}-V^{\text{s}}_k &\leq -\frac{\alpha_k}{L_\text{s}} \left(2-\alpha_k L_\text{s}\right)\left(1+\alpha_k L_\text{s}\right) \|\nabla f(x_{k+1})\|^2 \label{GD_LF_smooth_1}\\
F_k &\leq \frac{2 F_0 V^{\text{s}}_{0}}{2 V^{\text{s}}_{0}+ \frac{c_1}{L_\text{s}} F_0 k},  \label{GD_LF_smooth_3}\\
\|\nabla f(x_k)\| &\leq \frac{L_\text{s}}{c_2}\frac{\|x_0-x^\star\|}{k},  \label{GD_LF_smooth_2}
\end{align}
\hspace{0.01in}
\end{subequations}
\noindent where $c_1,c_2\in(0,\infty)$ are problem-independent constants.
\end{thm}
The proof builds on standard results \cite{Nesterov_Intro_2014,Teboulle2022}, but an analysis encompassing arbitrary sequences $(\alpha_k)_{k \in \N}$, and yielding (\ref{GD_LF_smooth}) is not present in the literature. Precisely: (\ref{GD_LF_smooth_1}) gives the existence of a Lyapunov function for (\ref{GD_method}) which can be used to show boundedness of the iterates and global convergence to the set $X^\star$;  (\ref{GD_LF_smooth_3}) and (\ref{GD_LF_smooth_2}) show convergence rates for function values and gradient which have no worse dependencies on $L_\text{s}$ and $k$ than those found in the literature and focusing on constant or smaller ranges of stepsizes \cite{Nesterov_Intro_2014,Teboulle2022}. 

The goal of this work is to develop a GD algorithm achieving similar guarantees to Theorem \ref{thm:GD_LF_smooth} under the standing Assumption \ref{Standing_ass} only. The only degree of freedom in (\ref{GD_method}) is the stepsize sequence, and we approach its design as an adaptive control problem where $(\alpha_k)_{k \in \N}$ is an input that can be chosen based on feedback and feedforward information to steer the iterate towards the set $X^\star$.

\section{Adapting stepsize to local smoothness}\label{sec:Adapting}

\subsection{A local smoothness estimate}\label{sec:Adapting-local}
The intuitive idea for using the GD method without assuming (global) smoothness of the objective is to adapt the stepsize sequence to the local geometry of the cost function. A natural measure of it along the GD iterates is the local smoothness estimate
\begin{equation}\label{def_L_k}
\begin{aligned}
L_k=L(x_{k+1},x_{k}) &\coloneqq \frac{\|\nabla f(x_{k+1})-\nabla f(x_{k})\|}{\|x_{k+1}-x_{k}\|}  , \\ 
&= \frac{\|\nabla f(x_k-\alpha_k \nabla f(x_k))-\nabla f(x_{k})\|}{\alpha_k\|\nabla f(x_{k})\|}. \\
\end{aligned}
\end{equation}
At iterate $k$, this estimate depends both on past information through $x_k$ (feedback) and one-step ahead future information through $x_{k+1}$ (feedforward). 

While after Theorem \ref{thm:GD_LF_smooth} it would be tempting to conjecture that a stepsize sequence satisfying $L_k \alpha_k \in (0,2)$ could satisfy our goal, the following results instill caution.
\newpage
\begin{lem}\label{lem:L_k_prop_1}
Consider a convex and differentiable $f$. 
\begin{enumerate}[(i)]
\item Given $L(y,x)$ defined in (\ref{def_L_k}), $\forall y,x \in \R^n$
\begin{equation}\label{L_k_prop_1}
f(y) \leq f(x)+\langle\nabla f(x),y-x\rangle + L(y,x)\|y-x\|^2.
\end{equation}
\item The sequence generated by (\ref{GD_method}) GD satisfies
\begin{equation}\label{L_k_prop_1_decrease}
F_{k+1} \leq F_k - \left(1-L_k \alpha_k \right) \alpha_k\|\nabla f(x_k)\|^2. 
\end{equation}
\item If $L_k \alpha_k \in \left(0,\frac{1}{2}\right)$, then
\begin{equation}\label{L_k_prop_1:LF}
V^{\text{s}}_{k+1}(x^\star)-V^{\text{s}}_k(x^\star) \leq - 2 \alpha_k F_{k+1}.
\end{equation}
\end{enumerate}
\end{lem}
Item (i) shows that, compared to the global smoothness constant $L_\text{s}$, the local smoothness estimate $L_k$ can be larger up to a factor of two, compare (\ref{f_glob_smooth_2}) and (\ref{L_k_prop_1}). A direct consequence of this is item (ii), where Eq. (\ref{L_k_prop_1_decrease}) clearly implies that a guaranteed function value decrease holds if $L_k \alpha_k \in (0,1)$. Finally, item (iii) shows that setting the stepsize to $L_k \alpha_k \in \left(0,\frac{1}{2}\right)$ guarantees the existence of the same Lyapunov function $V^{\text{s}}$ in (\ref{LF_smooth}).

We note that the sufficient condition $L_k \alpha_k \in (0,1)$ for function value decrease in item (ii) is, in general, also necessary. Indeed (\ref{L_k_prop_1_decrease}) follows immediately from (\ref{L_k_prop_1}), which has recently been shown to be tight \cite{Mishkin2024}.
\begin{lem}\label{lem:L_k_prop-2}\cite[Proposition 2.3]{Mishkin2024} 
Given $\beta \in(0,1)$, there exist a convex and differentiable $f_\beta$, $y,x \in \R^n$ such that
\begin{equation}\label{L_k_prop_2}
f_\beta(y) \geq f_\beta(x)+\langle\nabla f_\beta(x),y-x\rangle + \beta L(y,x)\|y-x\|^2. 
\end{equation}
\end{lem}

\subsection{Adaptive feedback-feedforward gradient descent}\label{sec:Adapting-AFFGD}
We propose here a novel stepsize update law that uses $L_k$ to adapt to the local geometry by combining feedback and feedforward mechanisms.
    
The update law reads as:
\begin{equation}\label{AFFGD_law}
\begin{aligned}
\alpha_k&= \min\biggl\{ \alpha_k^{(1)},\alpha_k^{(2)} \biggr\}, \\
 \text{where}\quad\alpha_k^{(1)}&=\frac{\gamma_k}{L_k},\quad \alpha_k^{(2)}=\frac{\alpha_{k-1}}{\gamma_k^2}\left(\frac{1-\gamma_k^2}{1-\gamma_{k-1}^2}\right),\\
(\gamma_k)_{k \in \N}&\subset (0,1).\\
\end{aligned}
\end{equation}
where $(\gamma_k)_{k \in \N}$ is a scalar sequence of parameters inside the specified range. 
The stepsize $\alpha_k$ is chosen as the smallest between the two upper bounds $\alpha_k^{(1)}$ and $\alpha_k^{(2)}$. The former one is the intuitive choice discussed in \ref{sec:Adapting-local}. It is worth observing that, whenever $\gamma_k \in \left(\frac{1}{2},1\right)$, $L_k \alpha_k^{(1)} \in \left(0,\gamma_k\right)$. That is, $\alpha_k$ in (\ref{AFFGD_law}) can be up to two times larger than the bound in item (iii) of Lemma \ref{lem:L_k_prop_1} guaranteeing the existence of $V^{\text{s}}$ in item (iii), and can become as large as the fundamental limit 
in item (ii) of Lemma \ref{lem:L_k_prop_1}.  
The bound $\alpha_k^{(2)}$ instead limits the increase of stepsize across two consecutive iterations and does not depend directly on the local geometry, but only on the last value of the stepsize. As it will be shown in Section \ref{sec:Adapting-robustness}, this bound also gives some inherent robustness to the algorithm. 
Intuitively, the first constraint is active in regions of the variable space where $f$ changes rapidly (or is \emph{less smooth}), whereas the second is active when, due to the function's \emph{flatness}, the stepsize would tend otherwise to overly increase. 
Finally, the parameters $(\gamma_k)_{k \in \N}$ are a tuning knob to navigate the speed of convergence vs. robustness trade-off discussed later. While any value (constant or time-varying) in $(0,1)$ is valid, one intriguing option is to use them as additional adaptive parameters. For example, they could be modified online so that the two upper bounds are as close as possible and thus $\alpha_k$ is maximized at every iteration. This option will be further explored in Section \ref{sec:Numerics}. 

The interconnection between the classic GD recursion (\ref{GD_method}) and the adaptation law (\ref{AFFGD_law}) is shown schematically in Figure \ref{block_diagram} and we will refer to in the following for brevity as AFFGD (adaptive feedback-feedforward gradient descent).
\begin{figure}[h!]
        \centering
        \resizebox{1\columnwidth}{!}{%
        
\input{./figures/block_diagram}
}
 \caption{The adaptive feedback-feedforward gradient descent algorithm (feedforward paths with dashed line).}
    \label{block_diagram}
    \end{figure}
    
The following result gives convergence guarantees for AFFGD by using a Lyapunov analysis which, albeit departing from the one used in Theorem \ref{thm:GD_LF_smooth} and under weaker requirements, yield qualitatively similar results.
\begin{thm}\label{thm:GD_LF_AFFGD}
Let $(x_k)_{k \in \N}$ be a sequence generated by the GD method (\ref{GD_method}) applied to a function $f$ satisfying Assumption \ref{Standing_ass} and with stepsize law (\ref{AFFGD_law}). Define the function
\begin{equation}\label{LF_local_smooth}
V^{\text{a}}_k(x^\star)\coloneqq\|x_k-x^\star\|^2+\frac{2 \alpha_{k-1}}{1-\gamma_{k-1}^2}F_k, \quad x^\star \in X^\star.
\end{equation}
Then for any $x^\star \in X^\star$ and $k \in \N$, it holds that
\begin{subequations}\label{GD_LF_local_smooth}
\begin{align}
V^{\text{a}}_{k+1}(x^\star)-V^{\text{a}}_k(x^\star) &\leq -\frac{2\gamma_{k}^2}{1-\gamma_{k}^2}\left(\alpha_k^{(2)}- \alpha_{k}\right)F_{k}-v_k \label{GD_LF_local_smooth_1}\\
F_k &\leq \frac{\|x_{0}-x^\star\|^2+2\alpha_{0}\frac{\gamma_0^2}{1-\gamma_0^2}F_0}{2\sum_{i=1}^{k-1}\alpha_i}\label{GD_LF_local_smooth_2}\\
\|\nabla f(x_k)\| &\xrightarrow{k\rightarrow \infty} 0  \label{GD_LF_local_smooth_3}
\end{align}
\end{subequations}
where
\begin{equation}\label{GD_LF_local_smooth_vk}
v_k:=\frac{\alpha_k}{L_{D_k}}\|\nabla f(x_k)\|^2+\frac{\alpha_k^2}{1-\gamma_k^2}\|\nabla f(x_{k+1})\|^2
\end{equation}
and $L_{D_k}$ is the local smoothness constant over a convex and compact set $D_k$ containing $x_{k}$ and $x^\star$.
\end{thm}
Eq. (\ref{GD_LF_local_smooth_1}) shows that function $V^{\text{a}}$ is a valid Lyapunov function for the closed-loop dynamics (\ref{GD_method})-(\ref{AFFGD_law}) which gives boundedness of the iterates, asymptotic optimality (\ref{GD_LF_local_smooth_3}) and global convergence to the set $X^\star$. Eq. (\ref{GD_LF_local_smooth_2}) gives a guaranteed last iterate convergence. As shown in the proof (cf. Eq. \ref{alpha_non_zero}), the stepsize sequence  $(\alpha_k)_{k \in \N}$ is separated from 0, and thus (\ref{GD_LF_local_smooth_2})  yields immediately a guaranteed convergence rate of $\mathcal{O}(\frac{1}{k})$, as in the standard smooth case. However, the denominator of (\ref{GD_LF_local_smooth_2}) points out that we can accelerate convergence by maximizing the sum of stepsizes. We can achieve this by adapting online the parameters $(\gamma_k)_{k \in \N}$ to make $\alpha_k^{(1)}$ and $\alpha_k^{(2)}$ as close as possible.
Compared to the recent literature on adaptive gradient descent \cite{Latafat2024_MathProg,Malitsky2024_AProx}, AFFGD provides convergence rate guarantees on the last iterate (\ref{GD_LF_local_smooth_2}), a Lyapunov function with the two standard terms relating to suboptimality distances (\ref{LF_local_smooth}), and a larger available upper bound on the stepsize with respect to the local geometry, compare with \cite[Table 1]{latafat24aL4DC}. It is also important on the other hand to recognize that the computation of $\alpha_k^{(1)}$ is a disadvantage of the proposed formulation as it involves forward prediction. While for some special cases this can be done without extra computation (e.g., when $f$ is quadratic with Hessian $M$, then $L_k$ only depends on $M$ and the current gradient), in general this requires a linesearch procedure that can be easily automated but might result in a more expensive per-iteration cost. Most importantly, we show next that limiting the growth rate of $\alpha_k$ (e.g., as currently done via $\alpha_k^{(2)}$) provides robustness to inexact gradients, for example due to errors in the linesearch procedure. 





\subsection{Robustness}\label{sec:Adapting-robustness}
The system theoretic view on AFFGD (Figure \ref{block_diagram}) prompts the question of robustness of the closed-loop. For example, one can consider stepsize updates where the first upper bound in (\ref{AFFGD_law}), involving forward prediction, is not exactly satisfied
\begin{equation}
\alpha_k\leq\tilde{\alpha}_k^{(1)}\coloneqq\frac{\gamma_k}{a_k L_k}, \quad a_k\in (0,1).
\end{equation}
Because $\tilde{\alpha}_k^{(1)}>\alpha_k^{(1)}$, this can capture errors in the gradient information (e.g., noisy evaluation, inexactness of the linesearch) inversely proportional to the parameter $a_k$. 

Let us define the scaled sequence $(\tilde{\gamma}_k)_{k \in \N}$ with \\$\tilde{\gamma}_k\coloneqq\frac{\gamma_k}{a_k}>\gamma_k$. 
If $a_k\in (\gamma_k,1)$, then $\tilde{\gamma}_k \in (\gamma_k,1) \forall k \in \N$. Then we can simply observe that we can still guarantee the results of Theorem \ref{thm:GD_LF_AFFGD} if we tighten the second upper bound in (\ref{AFFGD_law}) correspondingly, that is we impose
\begin{equation}
\tilde{\alpha}_k^{(2)}\coloneqq \frac{\alpha_{k-1}}{\tilde{\gamma}_k^2}\left(\frac{1-\tilde{\gamma}_k^2}{1-\tilde{\gamma}_{k-1}^2}\right).
\end{equation}
Indeed the conditions prescribed for the stepsize (\ref{AFFGD_law}) are satisfied with respect to the scaled sequence $(\tilde{\gamma}_k)_{k \in \N}\subset (0,1)$. This observation provides two insights. 
First, limiting the growth rate of $\alpha_k$ (through $\alpha_k^{(2)}$) adds robustness to inexact gradient information, which also contributes to understanding the role of the second upper bound (\ref{AFFGD_law}). 
Second, the tuning parameters $(\gamma_k)_{k \in \N}$ provide a means to navigate the trade-off between speed of convergence (when it is chosen adaptively to make $\alpha_k^{(1)}$ and $\alpha_k^{(2)}$ close and thus maximize the rate of convergence) and robustness (when it is chosen away from 1 to have robustness against perturbations $a_k\in (\gamma_k,1)$). 

It is natural to ask what happens when $a_k\in (0,\gamma_k]$, which models scenarios where perturbations are large or $\gamma_k$ is chosen close to 1.
In this case the analysis in Theorem \ref{thm:GD_LF_AFFGD} does not apply but the following result provides a first answer.
\begin{lem}\label{thm:GD_LF_AFFGD_pert}
Let $(x_k)_{k \in \N}$ be a sequence generated by the GD method (\ref{GD_method}) applied to a function $f$ satisfying Assumption \ref{Standing_ass} and with stepsize law 
\begin{equation}\label{GD_LF_local_smooth_pert_stepsize}
\alpha_k= \frac{\gamma_k}{a_k L_k},\quad \gamma_k \in (0,1), \; a_k\in (0,\gamma_k].
%
\end{equation}
Define the function
\begin{equation}\label{GD_LF_local_smooth_pert_candidate}
V^{\text{p}}_k(x^\star)\coloneqq\|x_k-x^\star\|^2+\frac{\alpha_{k-1}}{\alpha_{k}}\|x_{k+1}-x_k\|^2,\quad x^\star \in X^\star.  
\end{equation}
Then for any $x^\star \in X^\star$ and $k \in \N$, it holds that
\begin{equation}\label{GD_LF_local_smooth_pert}
V^{\text{p}}_{k+1}(x^\star)-V^{\text{p}}_k(x^\star)\leq -\left(\frac{\alpha_{k-1}^2}{\alpha_k^2} - \frac{\gamma_k^2}{a_k^2}\right)\|x_{k+1}-x_k\|^2-2\alpha_k F_{k+1}.
\end{equation}
\end{lem}
Note that Eq. (\ref{GD_LF_local_smooth_pert_stepsize}) allows perturbations to even determine stepsizes that results in $L_k \alpha_k>1$.
Condition (\ref{GD_LF_local_smooth_pert}) shows that, even in such extreme scenarios, limiting the growth rate of $\alpha_k$  guarantees \emph{robust} convergence. Indeed, if in addition to (\ref{GD_LF_local_smooth_pert_stepsize}) it holds 
\begin{equation}\label{GD_LF_local_smooth_pert_stepsize_2}
\alpha_k\leq \frac{a_k}{\gamma_k}\alpha_{k-1}
\end{equation}
then the decrease of $V^{\text{p}}$ is guaranteed. While (\ref{GD_LF_local_smooth_pert_stepsize_2}) is restrictive as for large perturbations it effectively prevents $\alpha_k$ from increasing, it provides another important characterization of the robustifying effect of limiting the growth rate even in this large perturbations regime. Moreover, we observe that requiring (\ref{GD_LF_local_smooth_pert_stepsize_2}) is not necessary because we are ignoring the second negative term on the r.h.s. of (\ref{GD_LF_local_smooth_pert}). We conjecture that positive growth rate conditions allowing $L_k \alpha_k>1$ are possible, but for space reasons we leave this for future work.


\section{Numerical Study}\label{sec:Numerics}

In this section we study numerically the performance of the proposed AFFGD algorithm and compare it with alternative methods from the literature. Codes to reproduce the results are available at the repository \footnote{\url{https://github.com/col-tasas/2025-AFFGD}}. 
We consider logistic regression, which is a convex optimization problem commonly used to train binary classifiers from data using large data sets.
Given $N$ features $s_i\in \R^n$ and labels $y_i=\pm1$, the goal is to find a linear classifier $x^\star \in \R^n$ by solving 
\begin{equation}\label{Log_ref}
	\min_{x \in \R^n} \quad\hspace*{-0.5em} \frac{1}{N}\sum_{1}^{N}\log(1+\exp(-y_ix^\top s_i)).
\end{equation}
The objective is globally smooth and the smallest smoothness constant is $L_\text{s}=\frac{1}{4N}\sigma_{\text{max}}(S)^2$, where $S\in \R^{N \times n}$ is the feature matrix. The solution of problem (\ref{Log_ref}) with GD has recently received attention \cite{Meng-arxiv24} due to the complex behavior of the iterates for \emph{large} stepsizes (i.e., greater than $\frac{2}{L_\text{s}}$) with not linearly-separable data. To reproduce this setting, we generate random data with $N=50$, $n=2$. 

\begin{figure}
\centering
    \includegraphics[width=1\columnwidth]{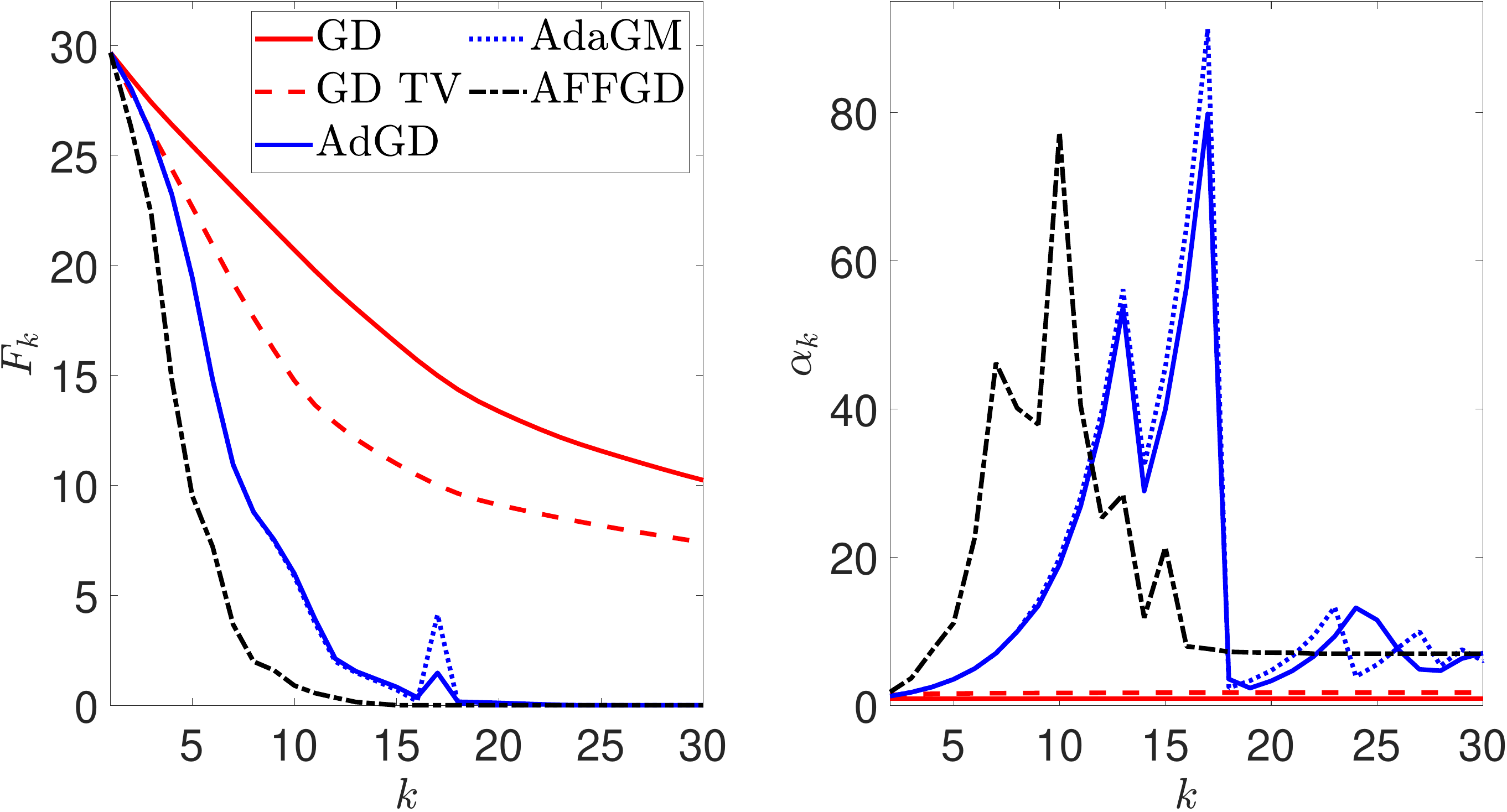}
\caption{Optimality gap and stepsizes of algorithms solving (\ref{Log_ref}).}
\label{Comparison_AGD}
\end{figure}

In a first set of results displayed in Figure \ref{Comparison_AGD}, we compare five GD algorithms (\ref{GD_method}) that differ for the step-size: \emph{GD} uses the classical choice $\alpha_k=\frac{1}{L_\text{s}}$ ($\simeq 1$ in this case); \emph{GD TV} is the dynamic update rule proposed in  \cite[Theorem 4]{Teboulle2022} whereby $L_\text{s}\alpha_k \in [1,2)$ and the stepisze is monotonically increased according to a pre-determined (or open-loop) law with $L_\text{s}\alpha_k \rightarrow 2$; \emph{AdGD} \cite[Algorithm 1]{Malitsky2024_AProx} and \emph{AdaGM} \cite[Algorithm 2]{Latafat2024_MathProg} are recently proposed adaptive GD schemes which also adapt to the local geometry by only using past information. Finally, \emph{AFFGD} is the update rule proposed in this work (\ref{AFFGD_law}) with a constant tuning parameter $\gamma=0.7$ and arbitrary initialization $\alpha_{-1}$.

The results show the clear impact that adaptation has on accelerating convergence. Even though the problem is globally smooth and the first two methods are guaranteed to asymptotically converge, they are slow compared to the adaptive ones. In that regard, the right plot shows that the sequence of adaptive stepsizes, during the iterations before convergence, sum up at a rate which is faster than linear. We also observe that \emph{AFFGD}, implemented here using a simple linesearch to verify the condition imposed by $\alpha_k^{(1)}$, markedly outperforms the other two methods at the cost of a slightly increased computational time (0.13$s$ against 0.1$s$). We finally note that, when simulating \emph{GD} with $\alpha>20\frac{1}{L_\text{s}}$ and \emph{GD TV} with $L_\text{s}\alpha_k \rightarrow 40$, we observed non-converging behaviours as described in \cite{Meng-arxiv24}. This is interesting because these stepsizes are still quite smaller than those (succesfully) employed in many iterations by the adaptive schemes (see the $y$ scale of the right plot in Figure \ref{Comparison_AGD}). This shows the importance of \emph{closed-loop} adaptation in optimization even for simple problems such as gradient descent applied to (\ref{Log_ref}).

\begin{figure}
\centering
    \includegraphics[width=1\columnwidth]{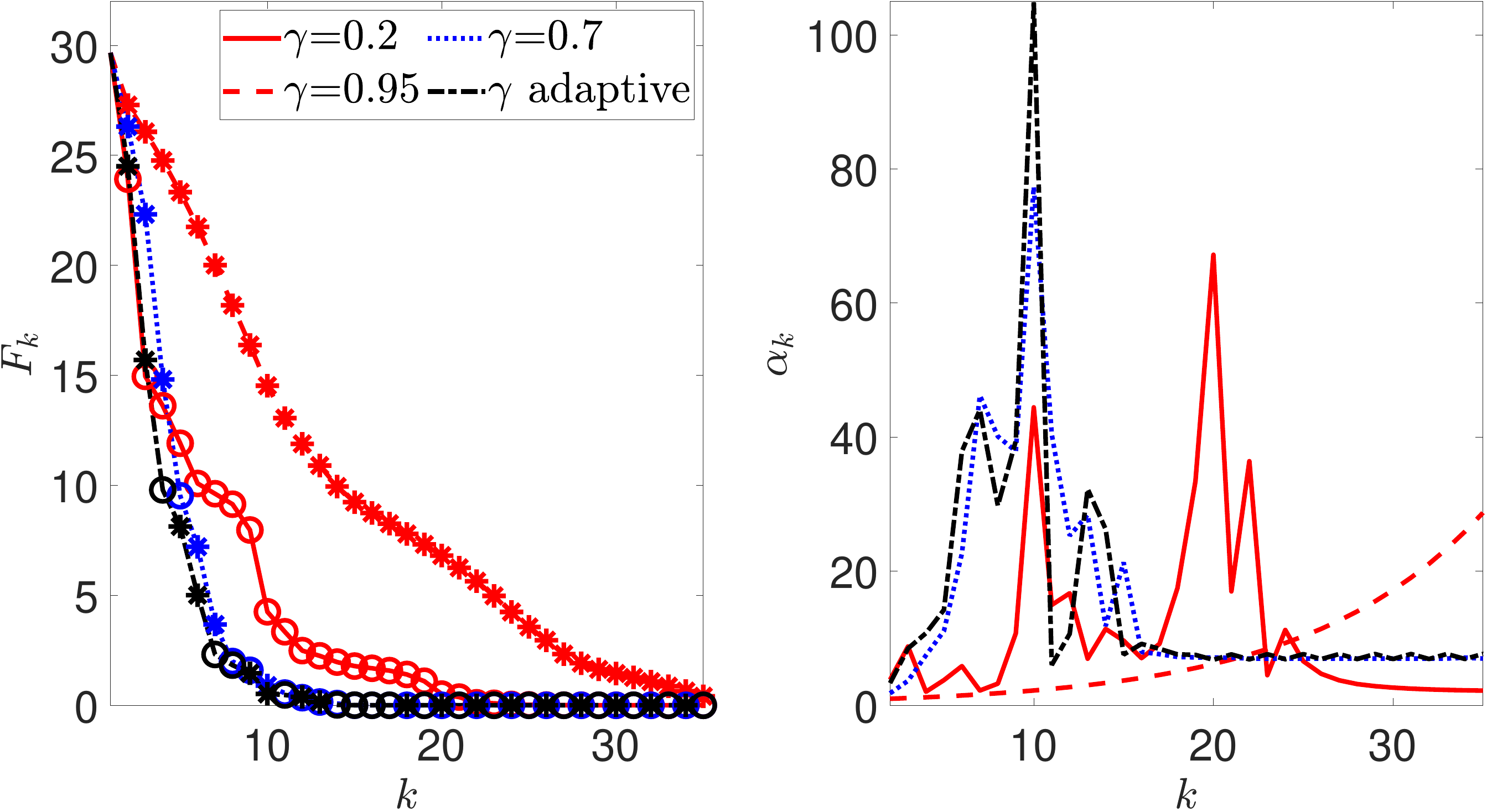}
\caption{Optimality gap and stepsizes of AFFGD for different choices of $\gamma_k$. Asterisk and circle markers denote points at which it holds that $\alpha_k=\alpha_k^{(1)}$ and $\alpha_k=\alpha_k^{(2)}$, respectively.}
\label{Comparison_AFFGD}
\end{figure}

We focus next in Figure \ref{Comparison_AFFGD} on AFFGD and investigate the effect of the sequence of tuning parameters $(\gamma_k)_{k \in \N}$. We compare three scenarios where this parameter is kept constant at some pre-defined value ($\gamma=\{0.2,0.7,0.95\}$) with the adaptive case where $\gamma_0=0.95$ and then it is changed adaptively using the simple recursion
\begin{equation}\label{eq:gamma_adaptive}
\begin{array}{rlll}
\gamma_k \!\!\!\!& = \begin{cases}
\frac{1}{\theta} \gamma_{k-1};
& \hspace{0.2in} \alpha_{k-1}=\alpha_{k-1}^{(1)}, \vspace{0.1in} \\ 
\theta\gamma_{k-1};
& \hspace{0.2in} \alpha_{k-1}=\alpha_{k-1}^{(2)}, \vspace{0.1in} \\ 
\end{cases},\quad k\in\mathbb{Z}_+
\end{array}
\end{equation}
where $\theta=0.9$ is a free parameter defining the strength of adaptation of $\gamma_k$. The rationale is to recursively update $\gamma_k$ based on the last active constraint in order to determine similar values for the two upper bounds $\alpha_k^{(1)}$ and $\alpha_k^{(2)}$, and by doing so maximize the sum of $\alpha_k$ which, as shown by our analyses (\ref{GD_LF_local_smooth_2}), accelerates the convergence rate.

The left plot shows, as expected,  that $\gamma=\{0.2,0.95\}$ yield low performance as they increase one of the bounds at the cost of strongly decreasing the other (which is then always active). 
On the contrary, the dynamic update rule (\ref{eq:gamma_adaptive}) is able to recover from the bad initialization $\gamma_0$ and determine larger values of stepsize and faster progress than those achieved with $\gamma=0.7$ (which was fine tuned offline).  

\begin{figure}[h]
\centering
    \includegraphics[width=1\columnwidth]{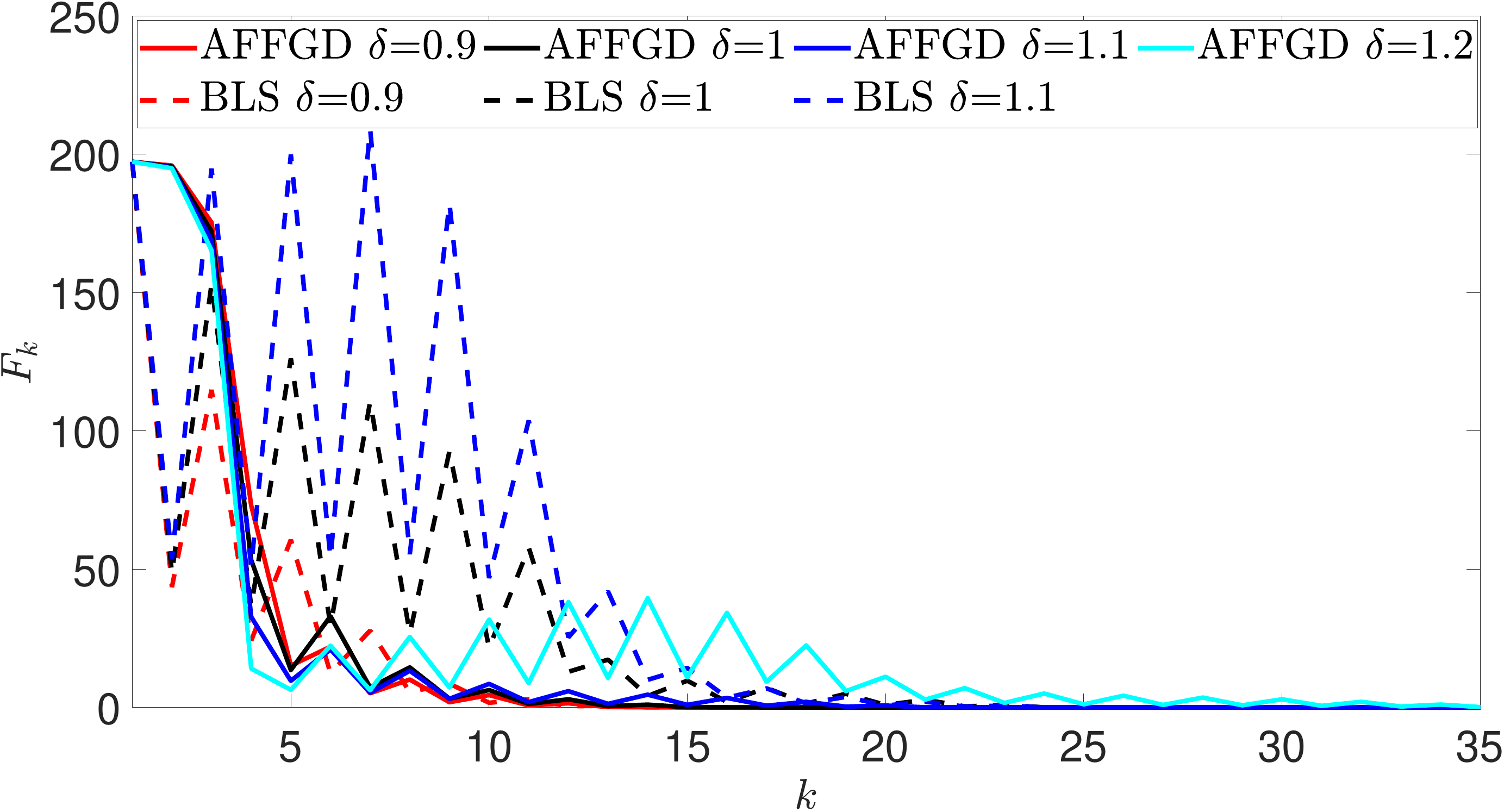}
\caption{Robustness of AFFGD vs. a simple backtracking strategy.}
\label{Comparison_AFFGD_BLS_robust}
\end{figure}

Finally, we take a numerical perspective on the robustness of AFFGD. Guided by the analytical results in Section \ref{sec:Adapting-robustness}, showing that limiting the growth rate of $\alpha_k$ through $\alpha_k^{(2)}$ robustifies the algorithm, we compare AFFGD with a simple backtracking line search (BLS) that sets $\alpha_k=\frac{\gamma_k}{L_k}$. We implement the update rule as $x_{k+1}=x_k-(1+\delta)\alpha_k \nabla f(x_k)$ to analyze the effect of numerical or gradient estimation errors quantified by the positive scalar $\delta$. Figure \ref{Comparison_AFFGD_BLS_robust} shows, in agreement with Lemma \ref{thm:GD_LF_AFFGD_pert}, that AFFGD (solid) is only marginally affected by such errors, while the convergence of BLS (dashed) degrades and is lost for $\delta>1.1$.
%
%
%
%

\section{CONCLUSIONS}

We consider the design of adaptive stepsize sequences for gradient descent methods driven by local properties of the objective. 
We frame the problem as an instance of adaptive controller where the input (stepsize) to the plant (GD method) is computed as a combination of feeedback and predicted information. 
Theoretical results and numerical experiments support the idea of pursuing closed-loop adaptation to accelerate convergence while increasing robustness.
Future directions include online optimization problems, where adaptation should capture both local geometry and instantaneous variations of the objective. 
We view the system theoretic framework in \cite{AI_Jakob_CDC_2025} as a promising starting point.

\section*{Appendix}
We report first an auxiliary result used later.
\begin{lem}\label{lem:f_glob_smooth}\cite[Theorem 5.8]{Beck_book_2017}
For any $L_\text{s}$-smooth and convex function $f$ it holds $\forall x,y \in \R^n$
\begin{equation}
f(y) \geq  f(x)+\langle\nabla f(x),y-x\rangle + \frac{1}{2L_\text{s}}\|\nabla f(y)-\nabla f(x)\|^2. 
\end{equation}
\end{lem}

\subsection{Proof of Theorem \ref{thm:GD_LF_smooth}}
\emph{Eq. (\ref{GD_LF_smooth_1}).}
For any $x^\star \in X^\star$ and any iteration $k\in \N$ of GD (\ref{GD_method}), apply the Pythagoras identity with $u=x_{k+1}-x^\star$ and $v=x_{k+1}-x_{k}$ 
\begin{equation}\label{thm:GD_LF_smooth_eq0}
\begin{aligned}
&V^{\text{s}}_{k+1}(x^\star)-V^{\text{s}}_k(x^\star)=\|x_{k+1}-x^\star\|^2-\|x_{k}-x^\star\|^2\\
&=2\langle x_{k+1}-x_{k},x_{k+1}-x^\star \rangle-\|x_{k+1}-x_{k}\|^2 \\
&=-2\alpha_k\langle \nabla f(x_{k}),x_{k}-x^\star \rangle+\alpha_k^2 \|\nabla f(x_{k})\|^2.\\
\end{aligned}
\end{equation}
We now bound the first term on the r.h.s using Lemma \ref{lem:f_glob_smooth} with $y=x^\star$ and $x=x_k$
\begin{equation}\label{thm:GD_LF_smooth_eq1}
-2\alpha_k\langle \nabla f(x_{k}),x_{k}-x^\star \rangle \leq -2\alpha_k F_k - \frac{\alpha_k}{L_\text{s}}\|\nabla f(x_k)\|^2, 
\end{equation}
which is a stronger bound than the one usually employed only using convexity. 
From $L_\text{s}$-smoothness of $f$ (\ref{f_glob_smooth_2}) with $y=x_{k+1}$ and $x=x_k$ we have
\begin{equation}\label{thm:GD_LF_smooth_eq1_1}
F_k \geq F_{k+1}+ (2-L_\text{s}\alpha_k)\frac{\alpha_k}{2}\|\nabla f(x_k)\|^2, 
\end{equation}
which, combined with (\ref{thm:GD_LF_smooth_eq1}) in (\ref{thm:GD_LF_smooth_eq0}) gives 
\begin{equation}\label{thm:GD_LF_smooth_eq2}
\begin{aligned}
&V^{\text{s}}_{k+1}-V^{\text{s}}_k\\
&\leq-2\alpha_k F_{k+1} -\frac{\alpha_k}{L_\text{s}} \left(1+\alpha_kL_\text{s}-(\alpha_kL_\text{s})^2\right) \|\nabla f(x_{k})\|^2.\\
\end{aligned}
\end{equation}
The second term on the r.h.s. is negative for $L_\text{s} \alpha_k \in (0,\frac{1+\sqrt{5}}{2})$. However, we can obtain a larger stepsize range by using again Lemma \ref{lem:f_glob_smooth} now with $y=x_{k+1}$ and $x=x^\star$ 
\begin{equation}\label{thm:GD_LF_smooth_eq2_2}
F_{k+1} \geq \frac{1}{2L_\text{s}}\|\nabla f(x_{k+1})\|^2, 
\end{equation}
and monotonicity of the gradient over the sequence of GD whenever $L_\text{s} \alpha_k \in (0,2)$  \cite[Lemma 2]{Teboulle2022}
\begin{equation}\label{thm:GD_LF_smooth_monot}
\|\nabla f(x_{k})\|^2 \geq \|\nabla f(x_{k+1})\|^2. 
\end{equation}
Using (\ref{thm:GD_LF_smooth_eq2_2})-(\ref{thm:GD_LF_smooth_monot}) in the r.h.s. of (\ref{thm:GD_LF_smooth_eq2}) gives
\begin{equation}
\begin{aligned}
V^{\text{s}}_{k+1}-V^{\text{s}}_k\leq -\frac{\alpha_k}{L_\text{s}} \left(2+\alpha_kL_\text{s}-(\alpha_kL_\text{s})^2\right) \|\nabla f(x_{k+1})\|^2,\\
\end{aligned}
\end{equation}
where the r.h.s. can be factorized as in (\ref{GD_LF_smooth_1}) and is negative for $L_\text{s} \alpha_k \in (0,2)$.

\emph{Eq. (\ref{GD_LF_smooth_3}).} Whereas (\ref{GD_LF_smooth_1}) does not lend itself to the classic telescopic argument used to show function value descrease, we can first observe that (\ref{GD_LF_smooth_1}) implies immediately that, for any $x^\star \in X^\star, x_0 \in \R^n$, there exists $d_0 \in [0,\infty)$ such that
\begin{equation}
\|x_{k}-x^\star\|^2\leq\|x_{0}-x^\star\|^2=:d_0, \quad k\in\N.
\end{equation}
Moreover, using again (\ref{thm:GD_LF_smooth_eq1}) and Cauchy-Schwartz we have  
\begin{equation}
F_k\leq d_0 \|\nabla f(x_k)\|. 
\end{equation}
Combining this with (\ref{thm:GD_LF_smooth_eq1_1}) yields 
\begin{equation}
F_{k+1} \leq F_{k} \left(1 - (2-L_\text{s}\alpha_k)\frac{\alpha_k}{2 d_0^2}F_{k}\right), 
\end{equation}
and thus
\begin{equation}\label{thm:GD_LF_smooth_fnc_value}
(F_{k+1})^{-1} - (F_{k})^{-1} \geq
 (2-L_\text{s}\alpha_k)\frac{\alpha_k}{2 d_0^2}, 
\end{equation}
where we also used $\frac{F_{k}}{F_{k+1}} \geq 1$ from the guaranteed objective decrease (\ref{thm:GD_LF_smooth_eq1_1}). We now sum (\ref{thm:GD_LF_smooth_fnc_value}) for $k=0,1,...,n-1$ and obtain 
\begin{equation}\label{thm:GD_LF_smooth_fnc_value_2}
(F_{n})^{-1} \geq(F_{0})^{-1} + \sum\limits_{k=0}^{n-1}(2-L_\text{s}\alpha_k)\frac{\alpha_k}{2 d_0^2}.
\end{equation}
We now make the following arguments. First, because  $L_\text{s} \alpha_k \in (0,2)$ and the technical restrictions on the limit behaviour, the second term on the r.h.s. can be lower bounded by a linear function of $n$ whose slope depends on the stepsize variations. Second (but not necessary for the qualitative result to hold), to capture the dependency of this term on $L_\text{s}$, we consider a linear dependence between $\alpha_k$ and $\frac{1}{L_\text{s}}$. This yields the existence of $c_1\in (0,\infty)$ such that 
\begin{equation}
(F_{n})^{-1} \geq(F_{0})^{-1}+ \frac{c_1}{2L_\text{s} d_0^2}n.
\end{equation}
Noting that $d_0=V^{\text{s}}_{0}$, this shows Eq. (\ref{GD_LF_smooth_3}).

\emph{Eq. (\ref{GD_LF_smooth_2}).}
We sum (\ref{GD_LF_smooth_1}) for $k=0,1,...,n-1$ and obtain
\begin{equation}
V^{\text{s}}_{0}-V^{\text{s}}_n\geq \sum\limits_{k=0}^{n-1} \frac{\alpha_k\left(2-\alpha_k L_\text{s}\right)\left(1+\alpha_k L_\text{s}\right)}{L_\text{s}} \|\nabla f(x_{k+1})\|^2.
\end{equation}
By using again the gradient monotonicity (\ref{thm:GD_LF_smooth_monot}) we have
\begin{equation}\label{thm:GD_LF_smooth_eq3}
\|\nabla f(x_{n})\|^2 \leq \frac{L_\text{s} \|x_{0}-x^\star\|^2}{\sum\limits_{k=0}^{n-1} \alpha_k\left(2-\alpha_k L_\text{s}\right)\left(1+\alpha_k L_\text{s}\right)}. 
\end{equation}
We now make 
similar arguments to those used for (\ref{thm:GD_LF_smooth_fnc_value_2}) to lower bound the denominator by a linear function of $n$ with linear dependence on $\frac{1}{L_\text{s}}$. It thus exists $c_2\in (0,\infty)$ such that 
\begin{equation}
\|\nabla f(x_{n})\|^2 \leq \frac{L_\text{s}^2 \|x_{0}-x^\star\|^2}{c_2^2 n} 
\end{equation}
and thus (\ref{GD_LF_smooth_2}) holds. 

\subsection{Proof of Lemma \ref{lem:L_k_prop_1}}
\emph{Eq. (\ref{L_k_prop_1}).} 
By convexity of $f$ (\ref{f_convex}), we have that $\forall x,y\in \R^n$
\begin{equation}
\begin{aligned}
f(y) &\leq f(x)+\langle\nabla f(y),y-x\rangle,\\
&= f(x)+\langle\nabla f(x),y-x\rangle+\langle \nabla f(y)- \nabla f(x),y-x \rangle,\\
& \leq f(x)+\langle\nabla f(x),y-x\rangle+ \| \nabla f(x)- \nabla f(y)\| \|y-x\|,\\
\end{aligned}
\end{equation}
where in the last inequality we used Cauchy-Schwarz. Plugging (\ref{def_L_k}) in the definition of $L(y,x)$ proves the statement.

\emph{Eq. (\ref{L_k_prop_1_decrease}).} 
This is a simple application of (\ref{L_k_prop_1}) taking $y=x_{k+1}$ and $x=x_k$ and using the GD method update
\begin{equation}
\begin{aligned}
f(x_{k+1}) &\leq f(x_k)+\langle\nabla f(x_k),x_{k+1}-x_k\rangle + L_k\|x_{k+1}-x_k\|^2,\\
 &\leq f(x_k)-(1- L_k \alpha_k)\alpha_k\|\nabla f(x_{k})\|^2,\\
\end{aligned}
\end{equation}
where $L_k=L(x_{k+1},x_{k})$.

\emph{Eq. (\ref{L_k_prop_1:LF}).}
Combine (\ref{thm:GD_LF_smooth_eq0}) with convexity of $f$ to get
\begin{equation}
V^{\text{s}}_{k+1}(x^\star)-V^{\text{s}}_k(x^\star)\leq -2 \alpha_k F_k+\alpha_k^2\|\nabla f(x_{k})\|^2.
\end{equation}
We upper bound the first term using (\ref{L_k_prop_1_decrease}) and get
\begin{equation}
V^{\text{s}}_{k+1}(x^\star)-V^{\text{s}}_k(x^\star)\leq -2 \alpha_k F_{k+1}+(2L_k\alpha_k-1)\alpha_k^2\|\nabla f(x_{k})\|^2
\end{equation}
which shows the sufficient decrease condition. 

\subsection{Proof of Theorem \ref{thm:GD_LF_AFFGD}}

\emph{Eq. (\ref{GD_LF_local_smooth_1}).} 
Like in the proof of Theorem \ref{thm:GD_LF_smooth}, we start off by applying Pythagoras identity (for the equality) and using Lemma \ref{lem:f_glob_smooth} (for the inequality)
\begin{equation}\label{thm:GD_LF_AFFGD_eq0}
\begin{aligned}
&\|x_{k+1}-x^\star\|^2-\|x_{k}-x^\star\|^2=-2\alpha_k\langle \nabla f(x_{k}),x_{k}-x^\star \rangle\\
&+\alpha_k^2 \|\nabla f(x_{k})\|^2 \leq -2\alpha_k F_k - \frac{\alpha_k}{L_{D_k}}\|\nabla f(x_k)\|^2+\alpha_k^2 \|\nabla f(x_{k})\|^2\\
\end{aligned}
\end{equation}
where $x^\star \in X^\star$. 
While we could have simply used convexity in the inequality (and drop the second term), local smoothness implies that Lemma \ref{lem:f_glob_smooth} holds $\forall x,y \in D_k$, 
where $D_k$ is a convex and compact set containing $x_{k}$ and $x^\star$, and $L_{D_k}$ is the associated local smoothness constant (\ref{f_loc_smooth}). We now work on the third term, where we use again Pythagoras identity but now with $u=x_{k+1}-x_{k}$ and $v=x_{k+2}-x_{k+1}$ yielding
\begin{equation}\label{thm:GD_LF_AFFGD_eq1}
\begin{aligned}
&\|x_{k+1}-x_{k}\|^2=\alpha_k^2 \|\nabla f(x_{k})\|^2=\alpha_k^2 \|\nabla f(x_{k+1})-\nabla f(x_{k})\|^2 \\
&-\alpha_k^2 \|\nabla f(x_{k+1})\|^2+2\alpha_k\langle \nabla f(x_{k+1}),x_{k}-x_{k+1} \rangle \\
&\leq \alpha_k^2 L_k^2 \|x_{k+1}-x_{k}\|^2-\alpha_k^2 \|\nabla f(x_{k+1})\|^2+2\alpha_k(F_k-F_{k+1})\\
\end{aligned}
\end{equation}
where for the inequality we used the definition of $L_k$ (first term) and convexity (third term).
Using now the upper bound on $\alpha_k$ due to $\alpha_k^{(1)}$ (\ref{AFFGD_law}), we get
\begin{equation}
\begin{aligned}
&(1-\gamma_k^2)\|x_{k+1}-x_{k}\|^2\leq-\alpha_k^2 \|\nabla f(x_{k+1})\|^2+2\alpha_k(F_k-F_{k+1}).\\
\end{aligned}
\end{equation}
Because  $(\gamma_k)_{k \in \N}\subset (0,1)$, we can divide the latter expression by $(1-\gamma_k^2)$ and plug this bound in (\ref{thm:GD_LF_AFFGD_eq0}) to obtain
\begin{equation}\label{thm:GD_LF_AFFGD_energy}
\begin{aligned}
&\|x_{k+1}-x^\star\|^2-\|x_{k}-x^\star\|^2\leq-\frac{ 2\alpha_k}{1-\gamma_k^2}F_{k+1}\\
&+\left(\frac{2\alpha_k}{1-\gamma_k^2}-2\alpha_k \right)F_{k}-v_k\\
\end{aligned}
\end{equation}
where $v_k$, defined in (\ref{GD_LF_local_smooth_vk}),
is a positive term for all $x_k \not\in X^\star$. Simple manipulations of (\ref{thm:GD_LF_AFFGD_energy}) yield
\begin{equation}
\begin{aligned}
&V^{\text{a}}_{k+1}(x^\star)-V^{\text{a}}_k(x^\star) \leq -2\left(\frac{\alpha_{k-1}}{1-\gamma_{k-1}^2} -\frac{\alpha_{k}\gamma_{k}^2}{1-\gamma_{k}^2}\right)F_{k}-v_k\\
&\overset{(\ref{AFFGD_law})}{=}-\frac{2\gamma_{k}^2}{1-\gamma_{k}^2}\left(\alpha_k^{(2)}- \alpha_{k}\right)F_{k}-v_k\\
\end{aligned}
\end{equation}
where the upper bound on $\alpha_k$ due to $\alpha_k^{(2)}$ (\ref{AFFGD_law}) guarantees negativity of the first term.

\emph{Eq. (\ref{GD_LF_local_smooth_2}).} 
We sum (\ref{thm:GD_LF_AFFGD_energy}) for $k=0,1,...,n-1$ and obtain
\begin{equation}\label{thm:GD_LF_AFFGD_energy_teles}
\begin{aligned}
&\|x_{n}-x^\star\|^2+\frac{2\alpha_{n-1}}{1-\gamma^2_{n-1}}F_n+\sum_{k=0}^{n-2}w_kF_{k+1} \\
&\leq \|x_{0}-x^\star\|^2+\left(\frac{2\alpha_{0}}{1-\gamma^2_{0}}-2\alpha_0 \right)F_0\\
\end{aligned}
\end{equation}
with
\begin{equation}
w_k:=2\left(\frac{\alpha_k}{1-\gamma_k^2}-\frac{\alpha_{k+1}}{1-\gamma_{k+1}^2}+\alpha_{k+1} \right).
\end{equation}
Observe now that, because of the upper bound on $\alpha_k$ due to $\alpha_k^{(1)}$ (\ref{AFFGD_law}) and $(\gamma_k)_{k \in \N}\subset (0,1)$, from  item (ii) of Lemma \ref{lem:L_k_prop_1} we have that for all iterations $F_k \geq F_{k+1}$. Note also that
\begin{equation}
\frac{2\alpha_{n-1}}{1-\gamma^2_{n-1}}+\sum_{k=0}^{n-2}w_k=\frac{2\alpha_{0}}{1-\gamma^2_{0}}+2\sum_{k=0}^{n-2}\alpha_{k+1}.
\end{equation}
Using these facts in (\ref{thm:GD_LF_AFFGD_energy_teles}) we finally obtain
\begin{equation}\label{thm:GD_LF_AFFGD_energy_final}
F_n\leq \frac{\|x_{0}-x^\star\|^2+2\alpha_{0}\frac{\gamma_0^2}{1-\gamma_0^2}F_0}{2\sum_{k=1}^{n-1}\alpha_k}.
\end{equation}
Note that the denominator of (\ref{thm:GD_LF_AFFGD_energy_final}) grows at least as fast as $k$ because we can show that the stepsize sequence  $(\alpha_k)_{k \in \N}$ is separated from 0. Indeed, the existence of the Lyapunov function $V^{\text{a}}$ for system (\ref{GD_method})-(\ref{AFFGD_law}) 
implies boundedness of the sequence $(x_k)_{k \in \N}$. For any $x^\star \in X^\star$, define $D^\star:=\overline{\text{conv}}(x^\star,x_0,x_1,..)$, which is closed and convex. Therefore, by local smoothness there exists bounded $L_{D^\star}$ such that (\ref{f_loc_smooth}) holds on $D^\star$. It is then 
\begin{equation}\label{alpha_non_zero}
\alpha_k \geq \frac{\gamma_k}{L_k}\geq \frac{\gamma_k}{L_{D^\star}}> 0, \quad \forall k \in \N
\end{equation}
which shows the desired property.

\emph{Eq. (\ref{GD_LF_local_smooth_3}).} 
We sum again (\ref{thm:GD_LF_AFFGD_energy}) for $k=0,1,...,n-1$ 
but this time also keep the last term $v_k$ (\ref{GD_LF_local_smooth_vk}). We obtain %
\begin{equation}\label{GD_LF_local_smooth_vk_2}
\sum_{k=1}^{n-1} \frac{\alpha_k}{L_{D_k}}\|\nabla f(x_k)\|^2 \leq c_3
\end{equation}
where the constant $c_3\in (0,\infty)$ exists due to the bounded quantities involved in (\ref{thm:GD_LF_AFFGD_energy}) and discussed in the previous item. Because $L_{D_k}\leq L_{D^\star}<\infty \forall k$, (\ref{GD_LF_local_smooth_vk_2}) gives summability of $\|\nabla f(x_k)\|^2$ and thus the result is proven. This shows that all cluster points of the sequences generated by (\ref{GD_method})-(\ref{AFFGD_law}) belong to $X^\star$. Using the Lyapunov function condition (\ref{GD_LF_local_smooth_1}) and classic fixed point arguments, one can then conclude that any sequence generated by the algorithm converges to a solution.

\subsection{Proof of Lemma \ref{thm:GD_LF_AFFGD_pert}}
Consider as starting point Eq. (\ref{thm:GD_LF_AFFGD_eq1}). Differently than in the proof of Theorem \ref{thm:GD_LF_AFFGD}, we can only use (\ref{GD_LF_local_smooth_pert_stepsize}) and thus we write in the last line  $\alpha_k^2 L_k^2 = \frac{\gamma_k^2}{a_k^2}$. Plugging this is in (\ref{thm:GD_LF_AFFGD_eq0}) and after simplifications we get
\begin{equation}\label{thm:GD_LF_AFFGD_pert_eq0}
\begin{aligned}
&\|x_{k+1}-x^\star\|^2-\|x_{k}-x^\star\|   \\
&\leq-2\alpha_k F_{k+1}- \alpha_k^2 \|\nabla f(x_{k+1})\|^2+\frac{\gamma_k^2}{a_k^2}\|x_{k+1}-x_{k}\|^2\\
&\leq-2\alpha_k F_{k+1}- \frac{\alpha_k^2}{\alpha_{k+1}^2} \|x_{k+2}-x_{k+1}\|^2+\frac{\gamma_k^2}{a_k^2}\|x_{k+1}-x_{k}\|^2.\\
\end{aligned}
\end{equation}
If we define the candidate Lyapunov function as in (\ref{GD_LF_local_smooth_pert_candidate}), simple manipulations of (\ref{thm:GD_LF_AFFGD_pert_eq0}) show that (\ref{GD_LF_local_smooth_pert}) holds.
\bibliographystyle{ieeetr}
\bibliography{biblio_AGD.bib} 

\end{document}

%% file: figures/block_diagram.tex
\begin{tikzpicture}
\node[draw, rectangle, minimum width=1.5cm, minimum height=1.5cm, text centered] at (1,2.5) (GD) {GD $(4)$};
\node[draw, rectangle, minimum width=1cm, minimum height=1cm, text centered] at (-0.3,0.5) (rk) {$\gamma_k$};
\node[draw, rectangle, minimum width=1cm, minimum height=0.6cm, text centered] at (-1.7,-0.6) (alpha2) {$\alpha^{(2)}$};
\node[draw, rectangle, minimum width=1cm, minimum height=0.6cm, text centered] at (-0.1,-1.5) (alpha1) {$\alpha^{(1)}$};
\node[draw, rectangle, minimum width=1.5cm, minimum height=1.5cm, text centered] at (1.7,-1.5) (fracL) {$\frac{1}{L_k}$};
\node[draw, rectangle, minimum width=1cm, minimum height=0.6cm, text centered] at (3.5,-1) (gradf) {$\nabla f$};
\node[draw, rectangle, minimum width=1cm, minimum height=0.6cm, text centered] at (3.5,-2) (gradfs) {$\nabla f^+$};
\node[draw, rectangle, minimum width=1.5cm, minimum height=1.5cm, text centered] at (-4.0,-1.05) (min) {$\min(\cdot,\cdot)$};
\node[draw, thick, dashdotted,rectangle, minimum width=10cm, minimum height=4cm] at(-0.2,-0.65) (red) {};
    \draw[->,thick] (min.west) -- ($(min.west)+(-1,0)$)|-  (GD.west);
    \draw[->,thick] ($(GD.east)$) --($(GD.east)+(3.5,0)$) |-($(4.8,-1)$)--(gradf.east);
    \draw[->,dashed,thick] ($(4.45,-1)$) |- (gradfs.east);
    \draw[->,thick] ($(gradf.west)$) -- ($(fracL.east)+(0,0.5)$);
    \draw[->,thick,dashed] ($(gradfs.west)$) -- ($(fracL.east)+(0,-0.5)$);
    \draw[->,thick] ($(fracL.west)$) -- (alpha1.east);
    \draw[->,thick] ($(rk.south)+(0.2,0)$) -- (alpha1.north);
    \draw[->,thick] ($(rk.south)-(0.2,0)$) |- (alpha2.east);
    \draw[->,thick] ($(alpha2.center)+(0,3.1)$) -- (alpha2.north);
    \draw[->,thick] (alpha2.west)--($(min.east)+(0,0.45)$);
    \draw[->,thick] (alpha1.west)--($(min.east)+(0,-0.45)$);
    \draw[->,thick] ($(alpha2.west)+(-0.3,0)$)|-($(rk.west)+(0,-0.2)$);
    \draw[->,thick] ($(alpha2.west)+(-0.7,-0.9)$)|-($(rk.west)+(0,0.2)$);
    \node at (2.3,2.8) {$x_k$};
    \node at (-5.5, -1.4) {$\alpha_k$};
\end{tikzpicture}